\documentclass{commat}

\allowdisplaybreaks
\usepackage{tikz-cd}
\DeclareMathOperator{\Hom}{Hom}

\title{%
    Cohomology, deformations and extensions of Rota-Baxter Leibniz algebras 
    }

\author{%
    Bibhash Mondal, Ripan Saha
    }

\affiliation{
    \address{Bibhash Mondal --
    Department of Mathematics, Behala College,
    Behala, 700060, Kolkata, India
        }
    \email{%
    mondaliiser@gmail.com
    }
    \address{Ripan Saha --
    Department of Mathematics, Raiganj University,  Raiganj 733134, West Bengal, India
        }
    \email{%
    ripanjumaths@gmail.com
    }
    
    }

\abstract{%
    A Rota-Baxter Leibniz algebra is a Leibniz algebra $(\mathfrak{g},[~,~]_{\mathfrak{g}})$ equipped with a Rota-Baxter operator $T : \mathfrak{g} \rightarrow \mathfrak{g}$. We define representation and dual representation of Rota-Baxter Leibniz algebras. Next, we define a cohomology theory of Rota-Baxter Leibniz algebras. We also study the infinitesimal and formal deformation theory of Rota-Baxter Leibniz algebras and show that  our cohomology is deformation cohomology. Moreover, We define an abelian extension of Rota-Baxter Leibniz algebras and show that equivalence classes of such extensions are related to the cohomology groups.
    }



\VOLUME{30}
\NUMBER{2}
\YEAR{2022}
\firstpage{93}
\DOI{https://doi.org/10.46298/cm.10295}

\begin{paper}

\section{Introduction}

A Rota-Baxter algebra is defined to be an associative algebra $A$ equipped with a linear map $T: A \to A$ satisfying Rota-Baxter identity: $T(x)T(y) = T(T(x)y + x T(y))$, for all $x, y\in A$. This Rota-Baxter identity is also known as Rota-Baxter equation of weight $0$, and $T$ is called the Rota-Baxter operator. Integration by parts is an example of a Rota–Baxter operator of weight $0$. Though it has a  natural connection with integral analysis, the Rota-Baxter algebra was not originated modeling this fact of  integral analysis, but was introduced by Glenn Baxter \cite{baxter} in his probability study of fluctuation theory. In a series of papers, Rota \cite{rota} derived some identities as well as gave some applications of Baxter algebra arising in probability and combinatorial theory. 
 
 A (left) Leibniz algebra \cite{loday93} is a vector space $\mathfrak{g}$ over a field F, equipped with a bilinear map  $[~,~]_{\mathfrak{g}} : \mathfrak{g}\otimes \mathfrak{g}\to \mathfrak{g}$ satisfying the (left) Leibniz identity:
   \[
   [x,[y,z]_{\mathfrak{g}}]_{\mathfrak{g}}
   = [[x,y]_{\mathfrak{g}},z]_{\mathfrak{g}} + [y,[x,z]_{\mathfrak{g}}]_{\mathfrak{g}},
   \quad \textup{for} \quad x,y,z \in \mathfrak{g}.
   \] 
   Leibniz algebras are often considered as non-commutative Lie algebras, since the Leibniz identity is equivalent to the Jacobi identity when the two-sided ideal $\lbrace x\in \mathfrak{g}~ \mid~ [x,x]=0\rbrace$ coincides with $\mathfrak{g}$. For this reason, a significant amount of research attempts to extend results on Lie algebras to Leibniz algebras. 
  
   Derivations on algebraic structures were first started by Ritt \cite{ritt} in the 1930s for commutative algebras and fields. The structure is called a differential (commutative) algebra. There is an enormous literature on this subject, including differential Galois theory, see \cite{M94}. In recent times, there is a numerous studies on different type of algebras with derivations, see \cite{das21}, \cite{guo saha}, \cite{TFS}. Similarly, degenerations of algebras is an interesting subject, which were studied in various papers, for the study of degenerations of Leibniz algebras and related structures, see \cite{KPPV}, \cite{IKV19}, \cite{KP19}. As Rota-Baxter operator is a kind of generalization of integral operator, therefore, it is natural to consider algebras with Rota-Baxter operator analogous to algebras with differentials. Jiang and Sheng \cite{JS} studied cohomology and abelian extensions of relative Rota-Baxter Lie algebras. In \cite{Guo}, the authors studied cohomology theory of Rota-Baxter Pre-Lie algebras of arbitrary weights. Tang, Sheng and Zhou \cite{TSZ} studied deformation theory of relative Rota-Baxter operator on Leibniz algebras. Recently, Das, Hazra, and Mishra \cite{Das} studied Rota-Baxter Lie algebras from cohomological point of view. 
   
    In this paper, our object of study is  Rota-Baxter Leibniz algebra. A Rota-Baxter Leibniz algebra is a Leibniz algebra $(\mathfrak{g},[~,~]_{\mathfrak{g}})$ equipped with a Rota-Baxter operator $T : \mathfrak{g} \rightarrow \mathfrak{g}$. We denote a Rota-Baxter Leibniz algebra by the notation $(\mathfrak{g}_T,[~,~]_{\mathfrak{g}})$. In this paper, we first define a representation of a Rota-Baxter Leibniz algebra. Given a Rota-Baxter Leibniz algebra $(\mathfrak{g}_T,[~,~]_{\mathfrak{g}})$, one can define a new bracket $[~,~]_\ast$ such that $(\mathfrak{g}_T,[~,~]_{\ast})$ is also Rota-Baxter Leibniz algebra. Now using this fact, if $(\mathfrak{g}_T,[~,~]_{\mathfrak{g}})$ be a Rota-Baxter Leibniz algebra and $(V,l_V,r_V,T_V)$ be a representation of it,  we get a new Rota-Baxter Leibniz algebra $(\mathfrak{g}_T,[~,~]_*)$ with representation $(V,l^{'}_V,r_V^{'},T_V)$ induced by the Rota-Baxter operator. We consider the Loday-Pirashvili cochain complex of this induced Leibniz algebra $(\mathfrak{g},[~,~]_*)$ with representation $(V,l^{'}_V,r_V^{'})$ to define the cohomology of Rota-Baxter Leibniz algebras. In the nth cochain group, we have two parts, one for Leibniz algebra and another one for the Rota-Baxter operator. In this paper, we also discuss dual representation a Rota-Baxter Leibniz algebra and nilpotent Rota-Baxter Leibniz algebra and obtain some interesting results. Next, we study one-parameter formal deformation theory of Rota-Baxter Leibniz algebra following  Gerstenhaber's \cite{G63}, \cite{G64} classical deformation theory for associative algebras. We study the the problem of extending a given deformation of order $n$ to a deformation of order $(n+1)$ and define the associated obstruction. We also study rigidity conditions for formal deformations. Finally, we define an abelian extension of a Rota-Baxter Leibniz algebra and show that how equivalence classes of such extensions are related to the cohomology groups.

This paper is organized as follows: In section \ref{sec1}, we recall Rota-Baxter operator, Leibniz algebra and its representation which we will use throughout the paper. In Section \ref{sec2}, we discuss some structural aspect of Rota-Baxter Leibniz algebras, also we discuss dual representation and nilpotent Rota-Baxter Leibniz algebras. In Section \ref{sec3}, we introduce the cohomology group of Rota-Baxter Leibniz algebra. In Section \ref{sec4}, we define Gerstenhaber's formal deformation theory of Rota-Baxter Leibniz algebras and showed that the cohomology defined in \ref{sec3} is a deformation cohomology. In the final Section, we discuss abelian extension and its relation with cohomology of Rota-Baxter Leibniz algebras.

\section{Preliminaries}\label{sec1}

 Let us recall the  following basic definitions from \cite{Das}, \cite{Loday}, \cite{Mandal}, \cite{Sheng}.

\begin{definition}
A \textbf{ Leibniz algebra} is a vector space $\mathfrak{g}$ together with a bilinear operation (called the bracket) $ [ ~,~]_{\mathfrak{g}}: \mathfrak{g} \times \mathfrak{g} \rightarrow \mathfrak{g}$ satisfying the following identity  \[ [x,[y,z]_{\mathfrak{g}}]_{\mathfrak{g}}=[[x,y]_{\mathfrak{g}},z]_{\mathfrak{g}}+[y,[x,z]_{\mathfrak{g}}]_{\mathfrak{g}} ,~~~~ \mbox{for}~ x,y,z \in \mathfrak{g}. \]
It is denoted by $(\mathfrak{g},[~,~]_{\mathfrak{g}}).$
\end{definition}
The above definition of Leibniz algebra is in fact the definition of left Leibniz algebra. In this paper, we will consider left Leibniz algebra simply as Leibniz algebra. Leibniz algebras are generalization of Lie-algebra. Any Lie algebra is a Leibniz algebra. A Leibniz algebra $(\mathfrak{g},[~,~])$ which satisfies $[a,a]=0 $ for all $a \in \mathfrak{g}$ is a Lie algebra.
\begin{example}
Let us consider the vector space $\mathbf{R}^2$ with standard basis $\{e_1,e_2\}$ with bracket defined by $[e_1,e_2]=0=[e_1,e_1],~ [e_2,e_1]=e_1=[e_2,e_2]$. Then $(\mathbf{R}^2,[~,~])$ is a Leibniz algebra.
\end{example}
\begin{definition}
Let $\mathfrak{g}=(g,[~,~]_{\mathfrak{g}})$ be a Leibniz algebra. A \textbf{Rota-Baxter} operator on $\mathfrak{g}$ is a linear map $ T:\mathfrak{g} \rightarrow \mathfrak{g}$ satisfying the following condition \[[T(x),T(y)]_{\mathfrak{g}}=T([T(x),y]_{\mathfrak{g}}+[x,T(y)]_{\mathfrak{g}}),~~~~\mbox{for all} ~~x, y \in \mathfrak{g}.\]
\end{definition}
\begin{definition}
A \textbf{Rota-Baxter Leibniz algebra} is a Leibniz algebra $(\mathfrak{g},[~,~]_{\mathfrak{g}})$ equipped with a Rota-Baxter operator $T : \mathfrak{g} \rightarrow \mathfrak{g}$. We denote a Rota-Baxter Leibniz algebra by the notation $(\mathfrak{g}_T,[~,~]_{\mathfrak{g}})$.
\end{definition}
\begin{example}
Consider the Leibniz algebra $(\mathbf{R}^2,[~,~])$ defined in example (2.2). Then the linear map $T: \mathbf{R}^2 \rightarrow \mathbf{R}^2,~ x \mapsto Ax$ where $A=\left(\begin{smallmatrix}
    0  &  b      \\
    0  &  0      
\end{smallmatrix}\right)$ for any $b \in \mathbf{R}$ is a Rota-Baxter operator on  $(\mathbf{R}^2,[~,~])$, hence $(\mathbf{R}^2_{T},[~,~])$ is a Rota-Baxter Leibniz algebra.
\end{example}

\begin{definition}
Let $(\mathfrak{g}_T,[~,~]_{\mathfrak{g}})$ and $(\mathfrak{g}_{T^{'}}^{'},[~,~]_{\mathfrak{g}^{'}})$ be two Rota-Baxter Leibniz algebras. A  \textbf{morphism} $\varphi : \mathfrak{g}_T \rightarrow \mathfrak{g}_{T^{'}}^{'}$ of Rota-Baxter Leibniz algebras is given by a Leibniz algebra homomorphism $\varphi : \mathfrak{g} \rightarrow \mathfrak{g}^{'}$ satisfying $T^{'}\circ \varphi= \varphi \circ T $. $\varphi$ is said to be an \textbf{isomorphism} if $\varphi$ is a linear isomorphism.
\end{definition}
Let $\mathfrak{g}_T$ be a Rota-Baxter Leibniz algebra. Then we denote Aut($\mathfrak{g}_T$) by the set of all isomorphism on $\mathfrak{g}_T$ .

\begin{definition}
Let $(\mathfrak{g},[~,~]_{\mathfrak{g}})$ be a Leibniz algebra. A $\mathbf{representation}$ of $(\mathfrak{g},[~,~]_{\mathfrak{g}})$ is a triple $(V,l_V,r_V)$ where $V$ is vector space together with  bilinear maps (called the left and right $\mathfrak{g}$-actions respectively ) $l_V : \mathfrak{g}\otimes V \rightarrow V ~~\mbox{and}~~ r_V : V\otimes \mathfrak{g} \rightarrow V $ satisfying the following conditions

\[l_V(x,l_V(y,u))=l_V([x,y]_{\mathfrak{g}},u)+l_V(y,l_V(x,u))\]
\[l_V(x,r_V(u,y))=r_V(l_V(x,u),y)+r_V(u,[x,y]_{\mathfrak{g}})\]
\[ r_V(u,[x,y]_{\mathfrak{g}})= r_V(r_V(u,x),y)+l_V(x,r_V(u,y))\]
for all $x,y \in \mathfrak{g}$ and $u\in V.$

\end{definition}
Now if we consider the vector space $V$ as $\mathfrak{g}$ itself with $l_{\mathfrak{g}}: \mathfrak{g} \times \mathfrak{g} \rightarrow \mathfrak{g}$, $r_{\mathfrak{g}}: \mathfrak{g} \times \mathfrak{g} \rightarrow \mathfrak{g}$ by $l_g(x,y)=[x,y]_{\mathfrak{g}}$ and $r_{\mathfrak{g}}(x,y)=[x,y]_{\mathfrak{g}}$  respectively for all $x,y \in \mathfrak{g}.$ Then $(\mathfrak{g},l_{\mathfrak{g}},r_{\mathfrak{g}})$ is a representation of $(\mathfrak{g},[~,~]_{\mathfrak{g}})$ which we call self representation.
Note that for self-representation the above three conditions reduce to the identity in the definition of Leibniz algebra.
\begin{definition}
Let $(\mathfrak{g}_T,[~,~]_{\mathfrak{g}})$ be a Rota-Baxter Leibniz algebra. A representation of $(\mathfrak{g}_T,[~,~]_{\mathfrak{g}})$ is a quadruple $(V,l_V,r_V,T_V) $, where $(V,l_V,r_V)$ is a representation of the Leibniz algebra $(\mathfrak{g},[~,~]_{\mathfrak{g}})$ and $T_V: V \rightarrow V$ is a linear map satisfying the following conditions
\[l_V(T(x),T_V(u))=T_V((l_V(T(x),u)+l_V(x,T_V(u)))\]
\[r_V(T_V(u),T(x))=T_V(r_V(T_V(u),x)+r_V(u,T(x)))\]
for all $x \in \mathfrak{g}$ and $u \in V.$
\end{definition}
Note that for a Rota-Baxter Leibniz algebra $(\mathfrak{g}_T,[~,~]_{\mathfrak{g}})$, the self representation of $(\mathfrak{g},[~,~]_{\mathfrak{g}})$ gives a representation $(\mathfrak{g},l_{\mathfrak{g}},r_{\mathfrak{g}},T)$ of the Rota-Baxter Leibniz algebra $(\mathfrak{g}_T,[~,~]_{\mathfrak{g}})$.

\section{Some Structural aspects of Rota-Baxter Leibniz algebras}\label{sec2}
\begin{proposition}
Let $(\mathfrak{g}_T,[~,~]_{\mathfrak{g}})$ be a  Rota-Baxter Leibniz algebra. Define
\[
[x,y]_*=[x,T(y)]_{\mathfrak{g}}+[T(x),y]_{\mathfrak{g}}
\quad \textup{for all} \quad x ,y  \in \mathfrak{g}.
\]
Then, 
\begin{enumerate}
\item  $(g,[~,~]_*)$ is a Leibniz algebra.
\item $T$ is also a Rota-Baxter operator on $(\mathfrak{g},[~,~]_*)$
\item The map $T:(g,[~,~]_{\mathfrak{g}}) \rightarrow (g,[~,~]_*)$ is a morphism of Rota-Baxter Leibniz algebra.
\end{enumerate}
\end{proposition}
 
\begin{proof}
\begin{enumerate}

\item 
It is easy to show that $[~,~]_*$ is a bilinear map.

Now, 
\begin{eqnarray*}
[x,[y,z]_{*}]_{*}
&=&[x,T([y,z]_{*})]_{\mathfrak{g}}+[T(x),[y,z]_{*}]_{\mathfrak{g}}\\
&=& [x,T([y,T(z)]_{\mathfrak{g}}+[T(y),z])]_{\mathfrak{g}}+[T(x),[y,T(z)]_{\mathfrak{g}}+[T(y),z]_{\mathfrak{g}}]_{\mathfrak{g}}\\
&=&[x,[T(y),T(z)]_{\mathfrak{g}}]_{\mathfrak{g}}+[T(x),[y,T(z)]_{\mathfrak{g}}]_{\mathfrak{g}}+[T(x),[T(y),z]_{\mathfrak{g}}]_{\mathfrak{g}}.
\end{eqnarray*}
Similarly, we have 
\[[[x,y]_{*},z]_{*}=[[x,T(y)]_{\mathfrak{g}},T(z)]_{\mathfrak{g}}+[[T(x),y]_{\mathfrak{g}},T(z)]_{\mathfrak{g}}+[[T(x),T(y)]_{\mathfrak{g}},z]_{\mathfrak{g}}\] and
\[[y,[x,z]_{*}]_{*}=[y,[T(x),T(z)]_{\mathfrak{g}}]_{\mathfrak{g}}+[T(y),[x,T(z)]_{\mathfrak{g}}]_{\mathfrak{g}}+[T(y),[T(x),z]_{\mathfrak{g}}]_{\mathfrak{g}}.\]
Next, using the identity in the definition of Leibniz algebra $(\mathfrak{g},[~,~]_{\mathfrak{g}})$ we have $[x,[y,z]_{*}]_{*}=[[x,y]_{*},z]_{*}+[y,[x,z]_{*}]_{*}$, for all $x,y,z \in \mathfrak{g}.$
\item  Now,
\begin{align*}
[T(x)&,T(y)]_{*} \\
&= [T(x),T(T(y))]_{\mathfrak{g}}+[T(T(x)),T(y)]_{\mathfrak{g}}\\
&= T([T(x),T(y)]_{\mathfrak{g}}+[x,T(T(y))]_{\mathfrak{g}})+T([T(T(x)),y]_{\mathfrak{g}}+[T(x),T(y)]_{\mathfrak{g}})\\
&= T([T(x),y]_*+[x,T(y)]_*).
\end{align*}
Therefore, $T$ is also a Rota-Baxter operator on the Leibniz algebra $(\mathfrak{g},[~,~]_*).$
\item  
It is obvious to observe that the map $T: (\mathfrak{g}_{T},[~,~]_{\mathfrak{g}}) \rightarrow (\mathfrak{g}_{T},[~,~]_*)$ , $g \mapsto T(g)$ is a morphism of Rota Baxter Leibniz algebra.
\qedhere
\end{enumerate} 
\end{proof}

\begin{proposition} 
Suppose $(V,l_{V},r_{V},T_V)$ be a representation of Rota-Baxter Leibniz algebra $(\mathfrak{g}_T,[~,~]_{\mathfrak{g}})$ and we define $l_V^{'} : \mathfrak{g} \otimes V \rightarrow V  ~, ~ r_{V}^{'} :  V  \otimes \mathfrak{g}  \rightarrow V $ respectively by 
\[l^{'}_V (x,u)=l_V(T(x),u)-T_V(l_V(x,u))\] 
\[r^{'}_V(u,x)=r_V(u,T(x))-T_V(r_V(u,x)),\]
 for all $x \in {\mathfrak{g}} , u \in V $. Then $(V,l_{V}^{'},r_V^{'},T_V)$ will be a representation of the Rota-Baxter Leibniz algebra $(\mathfrak{g}_T,[~,~]_*)$.
\end{proposition}
\begin{proof}
First, we will show that  $\left( V, l_V', r_V' \right)$ is in fact a representation of the Leibniz algebra $\left( \mathfrak{g}_T, [~,~]_* \right)$. For $x,y \in \mathfrak{g},u \in V$ we have,
\begin{align*}
l'_V(x,&l_V^{'}(y,u))-l_V^{'}([x,y]_*,u)-l_V^{'}(y,l_V^{'}(x,u))\\
={ }&{ }l_V(T(x),l_V^{'}(y,u))-(T_V\circ l_V)(x,l_V^{'}(y,u))-l_V(T([x,y]_*),u)+(T_V \circ l_V)([x,y]_*,u)\\
&-l_V(T(y),l_V^{'}(x,u))+(T_V\circ l_V)(y,l_V^{'}(x,u) \\
={ }&{ }l_V(T(x),l_V(T(y),u))-l_V(T(x),(T_V\circ l_V)(y,u))-(T_V\circ l_V)(x,l_V(T(y),u))\\
&+(T_V\circ l_V)(x,(T_V\circ l_V) (y,u))-l_V([T(x),T(y)]_{\mathfrak{g}},u)\\
&+(T_V\circ l_V)([x,Ty]_{\mathfrak{g}},u)+(T_V\circ l_V)([T(x),y]_{\mathfrak{g}},u)-l_V(T(y),l_V(T(x),u))\\
&+l_V(T(y),(T_V\circ l_V) (x,u))+(T_V\circ l_V)(y,l_V(T(x),u))-(T_V\circ l_V)(y,(T\circ l_V)(x,u))\\
={ }&{ }\bigg (l_V(T(x),l_V(T(y),u))-l_V([T(x),T(y)]_{\mathfrak{g}},u)-l_V(T(y),l_V(T(x),u)) \bigg )\\
&-l_V(T(x),(T_V\circ l_V)(y,u))\\
&-(T_V\circ l_V)(x,l_V(T(y),u))+(T_V\circ l_V)(x,(T_V\circ l_V) (y,u))+(T_V\circ l_V)([x,Ty]_{\mathfrak{g}},u)\\
&+(T_V\circ l_V)([T(x),y]_{\mathfrak{g}},u)\\
&+l_V(T(y),(T_V\circ l_V) (x,u))+(T_V\circ l_V)(y,l_V(T(x),u))-(T_V\circ l_V)(y,(T_V\circ l_V)(x,u))\\
={ }&{ }-T_V\bigg(l_V(T(x),l_V(y,u))+l_V(x,(T_V\circ l_V)(y,u))\bigg)-(T_V\circ l_V)(x,l_V(T(y),u))\\
&+(T_V\circ l_V)(x,(T_V\circ l_V) (y,u))+(T_V\circ l_V)([x,Ty]_{\mathfrak{g}},u)+(T_V\circ l_V)([T(x),y]_{\mathfrak{g}},u)\\
&+T_V\bigg(l_V(T(y),l_V(x,u))+l_V(y,(T_V\circ l_V)(x,u))\bigg)+(T_V\circ l_V)(y,l_V(T(x),u))\\
&-(T_V\circ l_V)(y,(T_V\circ l_V)(x,u))\\
={ }&{ }\bigg (-(T_V \circ l_V)(T(x),l_V(y,u))+(T_V\circ l_V)([T(x),y]_{\mathfrak{g}},u)+(T_V\circ l_V)(y,l_V(T(x),u))\bigg )\\
&\bigg (-(T_V\circ l_V)(x,l_V(T(y),u))+(T_V\circ l_V)([x,Ty]_{\mathfrak{g}},u)+(T_V \circ l_V)(T(y),l_V(x,u)) \bigg )\\
&+\bigg ( -(T_V \circ l_V)(x,(T_V\circ l_V)(y,u))+ (T_V \circ l_V)(x,(T_V\circ l_V)(y,u))\bigg )\\
&+\bigg ((T_V \circ l_V)(y,T_V \circ l_V(x,u))-(T_V \circ l_V)(y,T_V \circ l_V(x,u))\bigg )\\
={ }&{ }0.
\end{align*}

Therefore,
\[l_V^{'}(x,l_V^{'}(y,u))=l_V^{'}([x,y]_{*},u)+l_V^{'}(y,l_V^{'}(x,u))\]
holds.
Similarly it can be show that the following equations holds
\[l_V^{'}(x,r_V^{'}(u,y))=r_V^{'}(l_V^{'}(x,u),y)+r_V^{'}(u,[x,y]_{*})\]
\[ r_V^{'}(u,[x,y]_{*})= r_V^{'}(r_V^{'}(u,x),y)+l_V^{'}(x,r_V^{'}(u,y))\]

Now,
\begin{align*}
T_V(l_V^{'}&(T(x),u)+l_V^{'}(x,T_V(u)))\\
={ }&{ }T_V(l_V(T(T(x)),u))-T_V(l_V(T(x),u))+l_V(T(x),T_V(u))-T_V(l_V(x,T_V(u)))\\
={ }&{ }T_V(l_V(T(T(x)),u))-T_V(l_V(T(x),u))+T_V(l_V(T(x),u))+T_V(l_V(x,T_V(u)))\\
&-T_V(l_V(x,T_V(u)))\\
={ }&{ }T_V(l_V(T(T(x)),u))\\
={ }&{ }T_V(l_V(T(T(x)),u))+T_V(l_V(T(x),T_V(u)))-T_V(l_V(T(x),T_V(u)))\\
={ }&{ }l_V(T(T(x)),T_V(u))-T_V(l_V(T(x),T_V(u)))\\
={ }&{ }l_V^{'}(T(x),T_V(u)).
\end{align*}
Therefore, 
\[ l_V^{'}(T(x),T_V(u))=T_V(l_V^{'}(T(x),u)+l_V^{'}(x,T_V(u))).\]
Similarly, it can be shown that 
\[r_V^{'}(T_V(u),T(x))=T_V(r_V^{'}(T_V(u),x)+r_V^{'}(u,T(x)))\]
for all $x \in \mathfrak{g}$ and $u \in V.$
Hence, $(V,l_V^{'},r_V^{'},T_V )$ is a representation of the Rota-Baxter Leibniz algebra $(\mathfrak{g}_T,[~,~]_*).$
\end{proof}
\subsection{Dual representation of a Rota-Baxter Leibniz algebras}\label{sec21}
Let $(V,l_V,r_V)$ be a representation of a Leibniz algebra $(\mathfrak{g},[~,~]_{\mathfrak{g}})$. Then we can get a representation $(V^*,l_{V^*},r_{V^*})$ of $(\mathfrak{g},[~,~]_{\mathfrak{g}})$, where $V^*$ is the dual vector space of $V$ and $l_{V^*}:  \mathfrak{g}\otimes V^* \rightarrow V^*,~ r_{V^*}: V^{*} \otimes \mathfrak{g} \rightarrow V^*$ defined by 
\[l_{V^*}(x,f_V)(u)=-f_V(l_V(x,u)) ~~~ \mbox{and}~~~ r_{V^*}(f_V,x)(u)=f_V(l_V(x,u)+r_V(u,x)),\]
for all $x\in \mathfrak{g}, f_V \in V^*, u \in V.$ For details see \cite{Sheng}. We can use this notion to get dual representation for Rota-Baxter Leibniz algebra.

\begin{proposition}
Let $(V,l_V,r_V,T_V)$ be a representation of a Rota-Baxter Leibniz algebra $(\mathfrak{g}_T,[~,~]_{\mathfrak{g}}).$ Now define 
\[l_{V^*}(x,f_V)(u)=-f_V(l_V(x,u)) ~~~ \mbox{and}~~~ r_{V^*}(f_V,x)(u)=f_V(l_V(x,u)+r_V(u,x))\]
for all $x\in \mathfrak{g}, f_V \in V^*, u \in V.$ Then $(V^*,l_{V^*},r_{V^*},-T_{V}^*)$ is a representation of the Rota-Baxter Leibniz algebra $(\mathfrak{g}_T,[~,~]_{\mathfrak{g}}).$ We call it the dual representation of Rota-Baxter Leibniz algebra $(\mathfrak{g}_T,[~,~]_{\mathfrak{g}}).$
\end{proposition}
\begin{proof}
Now for all $x\in \mathfrak{g}, f_V \in V^*, u \in V.$
\begin{align*}
&l_V{^*}(x,l_V{^*}(y,f_V))(u)-l_V{^*}([x,y]_{\mathfrak{g}},f_V)(u)-l_V{^*}(y,l_V{^*}(x,f_V))(u)\\
&=-l_V{^*}(y,f_V)(l_V(x,u))+ f_V(l_V([x,y]_{\mathfrak{g}},u))+l_V{^*}(x,f_V)(l_V(y,u))\\
&=f_V(l_V(y,l_V(x,u)))+f_V(l_V([x,y]_{\mathfrak{g}},u))-f_V(l_V(x,l_V(y,u)))\\
&=-f_V\bigg(l_V(x,l_V(y,u))-l_V([x,y]_{\mathfrak{g}},u)-l_V(y,l_V(x,u)) \bigg) \\
&=0.
\end{align*}
Similarly, it can be shown that
\[l_{V^*}(x,r_{V^*}(u,y))=r_{V^*}(l_{V^*}(x,u),y)+r_{V^*}(u,[x,y]_{\mathfrak{g}})\]
\[ r_{V^*}(u,[x,y]_{\mathfrak{g}})= r_{V^*}(r_{V^*}(u,x),y)+l_{V^*}(x,r_{V^*}(u,y)).\]

Again,
\begin{align*}
&l_{V^*}(T(x),-T^*_V(f_V))(u)+T_V^*\bigg(l_{V^*}(T(x),f_V)+l_{V^*}(x,-T_V^*(f_V))\bigg)(u)\\
&=T_V^*(f_V)(l_V(T(x),u))+l_V^*(T(x),f_V)(T_V(u))+l_{V^*}(x,-T_V^*(f_V))(T_V(u))\\
&=f_V(T_V(l_V(T(x),u)))-f_V(l_V(T(x),T_V(u)))+T_V^*(f_V)(l_V(x,T_V(u)))\\
&=f_V(T_V(l_V(T(x),u)))-f_V(l_V(T(x),T_V(u)))+f_V(T_V(l_V(x,T_V(u))))\\
&=f_V\bigg(T_V\bigg(l_V(T(x),u)+l_V(x,T_V(u))\bigg)-l_V(T(x),T_V(u))\bigg)\\
&=0.
\end{align*}
Therefore, we have $l_{V^*}(T(x),-T^*_V(f_V))=-T_V^*\left(l_{V^*}(T(x),f_V)+l_{V^*}(x,-T_V^*(f_V))\right)$, for all $x\in \mathfrak{g}, f_V \in V^*$. Now, 
\begin{align*}
r_{V^*}(&-T_V^*(f_V),T(x))(u)+T_V^*\bigg(r_{V^*}(f_V,T(x))+r_{V^*}(-T_V^*(f_V),x)\bigg)(u)\\
={ }&{ }-T_V^*(f_V)\bigg(l_V(T(x),u))+r_V(u,T(x))\bigg)+r_{V^*}(f_V,T(x))(T_V(u))\\
&+r_{V^*}(-T_V^*(f_V),x)(T_V(u))\\
={ }&{ }-f_V\bigg(T_V\big(l_V(T(x),u)+r_V(u,T(x))\big)\bigg)+f_V\bigg(r_V(T_V(u),T(x))+l_V(T(x),T_V(u))\bigg)\\
&-T_V^*(f_V)\bigg(l_V(x,T_V(u))+r_V(T_V(u),x)\bigg)\\
={ }&{ }-f_V\bigg(T_V\big(l_V(T(x),u)+r_V(u,T(x))\big)\bigg)+f_V\bigg(r_V(T_V(u),T(x))+l_V(T(x),T_V(u))\bigg)\\
&-f_V\bigg(T_V\big(l_V(x,T_V(u))+r_V(T_V(u),x)\big)\bigg)\\
={ }&{ }f_V\bigg(r_V(T_V(u),T(x))-T_V\bigg(r_V(T_V(u),x))+r_V(u,T(x))\bigg)\bigg) \\
&+f_V\bigg(l_V(T(x),T_V(u)) - T_V\bigg(l_V(T(x),u)+l_V(x,T_V(u))\bigg)\bigg)\\
={ }&{ }0.
\end{align*}
Therefore, $r_{V^*}(-T^*(f_V),T(x))=-T_V^*\left(r_{V^*}(f_V,T(x))+r_{V^*}(-T_V^*(f_V),x)\right)$ holds.
\end{proof}
\subsection{Nilpotent operator on Rota-Baxter Leibniz algebra}\label{sec22}
\begin{definition}
A Rota-Baxter operator $T: \mathfrak{g} \to \mathfrak{g}$ is called nilpotent if there exist a positive integer $n$ such that $T^n=0$. The smallest such $n$ is called the degree of nilpotency.
\end{definition}
\begin{example}
Let us consider the three dimensional vector space $\mathbf{R}^3$ with standard basis $\{e_1,e_2,e_3\}$ and the  bracket defined by $[e_3,e_2]=e_2, [e_3,e_1]=e_1+e_2$. Then $(\mathbf{R}^3,[~,~])$ is a Leibniz algebra. Now define  a linear map $T: \mathbf{R}^3 \rightarrow \mathbf{R}^3,~ x \mapsto Ax$ where $A=\left(\begin{smallmatrix}
    0 &0 &  b      \\
    0  & 0 & c \\
    0& 0 &0      
\end{smallmatrix}\right)$ for any $b ,c\in \mathbf{R}$. Then $T$ is a nilpotent Rota-Baxter operator on the Rota-Baxter Leibniz algebra  $(\mathbf{R}^3_{T},[~,~])$. 
\end{example}

\begin{definition}
Let $(\mathfrak{g}_T,[~,~]_{\mathfrak{g}})$ be a  Rota-Baxter Leibniz algebra. Define 
\begin{gather*}
[x,y]_0=[x,y]_{\mathfrak{g}}, \qquad
[x,y]_1=[x,y]_{*}=[x,Ty]_0+[Tx,y]_0,\\
[x,y]_r=[x,Ty]_{r-1}+[Tx,y]_{r-1},
\end{gather*}
for all $r\in \mathbb{N}$ and $x,y  \in \mathfrak{g}$.
\end{definition}
\begin{proposition}
Let $(\mathfrak{g}_T,[~,~]_{\mathfrak{g}})$ be a  Rota-Baxter Leibniz algebra. Then
\[
[x,y]_n= \sum_{r=0}^n \binom{n}{r}[T^{n-r}(x),T^r(y)]_{\mathfrak{g}}
\]
for all $n \in \mathbb{N}$, where $T^{0}=Id_{\mathfrak{g}}$.
\end{proposition}
\begin{proof}
We will prove this result by Mathematical Induction on $n.$

For $n=1$ we have, $[x,y]_1= {\binom{1}{0} }[Tx,y]_{\mathfrak{g}}+{\binom{1}{1} }[x,Ty]_{\mathfrak{g}}=[Tx,y]_{\mathfrak{g}}+[x,Ty]_{\mathfrak{g}}$ which is true by definition. We assume the statement is true for $n=m$, $m\in \mathbb{N}$ , therefore,
\begin{align*}
[x,y&]_{m+1}\\
={ }&{ }[Tx,y]_m+[x,Ty]_m \\
={ }&{ }\sum_{r=0}^m{\binom{m}{r}}[T^{m-r}(Tx),T^r(y)]_{\mathfrak{g}}+\sum_{r=0}^m{\binom{m}{r}}[T^{m-r}(x),T^r(Ty)]_{\mathfrak{g}}\\
={ }&{ }{\binom{m}{0}}[T^{m+1}(x),y]_{\mathfrak{g}}+\sum_{r=1}^m{\binom{m}{r}}[T^{m-r}(Tx),T^r(y)]_{\mathfrak{g}}+\sum_{r=0}^{m-1}{\binom{m}{r}}[T^{m-r}(x),T^r(Ty)]_{\mathfrak{g}}\\
&+{\binom{m}{m}}[x,T^{m+1}(y)]_{\mathfrak{g}} \\
={ }&{ }[T^{m+1}(x),y]_{\mathfrak{g}}+\sum_{r=1}^m{\binom{m}{r}}[T^{m+1-r}(x),T^r(y)]_{\mathfrak{g}}+\sum_{r=1}^{m}{\binom{m}{r}-1}[T^{m+1-r}(x),T^r(y)]_{\mathfrak{g}}\\
&+[x,T^{m+1}(y)]_{\mathfrak{g}}\\
={ }&{ }[T^{m+1}(x),y]_{\mathfrak{g}}+\sum_{r=1}^m{m+\binom{1}{r}}[T^{m+1-r}(x),T^r(y)]_{\mathfrak{g}}+[x,T^{m+1}(y)]_{\mathfrak{g}}\\
={ }&{ }\sum_{r=0}^{m+1}{m+\binom{1}{r}}[T^{m+1-r}(x),T^r(y)]_{\mathfrak{g}}.
\qedhere
\end{align*}
\end{proof}
\begin{corollary}
If $T$ is idempotent, that is, $T^2=T$, then
\[
[x,y]_{\mathfrak{n}}=[x,y]_1+(2^n-2)[Tx,Ty]_{\mathfrak{g}}
\quad \textup{for all} \quad n\in \mathbb{N}.
\]
\end{corollary}
\begin{proof}
Observe that for all $n\in \mathbb{N}$, we have
\begin{align*}
[x,y]_{\mathfrak{n}}&=\sum_{r=0}^n{\binom{n}{r}}[T^{n-r}(x),T^r(y)]_{\mathfrak{g}}\\
&=[Tx,y]_{\mathfrak{g}}+\{{\binom{n}{1}}+{\binom{n}{2}}+{\binom{n}{3}}+\dotsb+{\binom{n}{n}-1}\}[Tx,Ty]_{\mathfrak{g}}+[x,Ty]_{\mathfrak{g}}\\
&=[x,y]_1+(2^n-2)[Tx,Ty]_{\mathfrak{g}}.
\qedhere
\end{align*}
 \end{proof}
\begin{proposition}
Let $(\mathfrak{g}_T,[~,~]_{\mathfrak{g}})$ be a  Rota-Baxter Leibniz algebra and $T : \mathfrak{g} \rightarrow \mathfrak{g}$ be an injective homomorphism. Then $[x,y]_n=[x,y]_{\mathfrak{g}}$ for all $n \in \mathbb{N},~ x,y \in \mathfrak{g}.$  
\end{proposition}
 
\begin{proof}
Now, $T([x,y]_{\mathfrak{g}})=[T(x),T(y)]_{\mathfrak{g}}=T([x,Ty]_0+[Tx,y]_0)=T([x,y]_{1}) $. Therefore,  
$[x,y]_1=[x,y]_{\mathfrak{g}}.$ By induction, we obtain that $[x,y]_n=[x,y]_{\mathfrak{g}}$ for all $n \in \mathbb{N}$, $x,y \in \mathfrak{g}.$   
\end{proof}
\begin{proposition}
Let $(\mathfrak{g}_T,[~,~]_{\mathfrak{g}})$ be a  Rota-Baxter Leibniz algebra and $T : \mathfrak{g} \rightarrow \mathfrak{g}$ is an nilpotent operator of index $n$. Then $[x,y]_{k}=0$ for any $k\geq (2n+1),~ x,y \in \mathfrak{g}.$ 
\end{proposition}
\begin{proof}
Now, $[x,y]_{2n+1}=\sum_{r=0}^{2n+1}{2n+\binom{1}{r}}[T^{2n+1-r}(x),T^r(y)]_{\mathfrak{g}}$, Now since each term of this sum contains some bracket $[T^{l}(x),T^{m}(y)]$ such that either $l$ or $m$ is greater than or equal to $n$. Hence $[x,y]_{2n+1}=0$. Similar arguments holds for any $k> (2n+1).$
\end{proof}
\begin{proposition}
Let $(\mathfrak{g}_T,[~,~]_{\mathfrak{g}})$ be a  Rota-Baxter Leibniz algebra and $T : \mathfrak{g} \rightarrow \mathfrak{g}$ be a surjective nilpotent operator. Then $[x,y]_{n}=0$ for all $n \in \mathbb{N},~ x,y \in \mathfrak{g}.$
\end{proposition}
\begin{proof}
Since $T$ is a nilpotent operator, there exists $m \in\mathbb{N}$ such that $T^m(x)=0$ for all $ x\in \mathfrak{g}$. Now, for any $x,y \in \mathfrak{g},~[x,y]_1=[Tx,y]_{\mathfrak{g}}+[x,Ty]_{\mathfrak{g}}$. Since $T$ is surjective hence there exists $x_1,y_1 \in \mathfrak{g}$ such that $Tx_1=x,Ty_1=y$. Repeatedly using the surjectivity of $T$ we get some  $x^{'},y^{'} \in \mathfrak{g}$ such that $T^{m-1}(x^{'})=x,T^{m-1}(y^{'})=y$. Then
\[
[x,y]_1
= [Tx,y]_{\mathfrak{g}} + [x,Ty]_{\mathfrak{g}}
= [T^{m}(x),y]_{\mathfrak{g}} + [x,T^{m}(y)]_{\mathfrak{g}}=0.
\]
Using Mathematical Induction we get, $[x,y]_n=0$ for all $n \in \mathbb{N},~~x,y \in \mathfrak{g}.$
\end{proof}

\section{Cohomology of Rota-Baxter Leibniz algebra}\label{sec3}
Let $(\mathfrak{g},[~,~]_{\mathfrak{g}})$ be a Leibniz algebra and $(V,l_V,r_V)$ be a representation of it.  For each  $n \geq 0$,  define $C^n_{LA}(\mathfrak{g},V)$ to be the  abelian group $\Hom(\mathfrak{g}^{\otimes n},V)$ and $\delta ^n$ to be the map $\delta^n : C^n_{LA}(\mathfrak{g},V) \to C^{n+1}_{LA}(\mathfrak{g},V)$ given by
\begin{align*}
 (\delta^n(f))&(x_1,x_2,\ldots ,x_{n+1})\\
 ={ }&{ }\sum_{i=1}^{n}(-1)^{i+1}l_{V}(x_i,f(x_1,\ldots,\hat{x_i},\ldots, x_{n+1}))+(-1)^{n+1}r_{V}(f(x_1,\ldots,x_n),x_{n+1}) \\
 &+ \sum_{1 \leq i <j \leq n+1}(-1)^if(x_1,\ldots, \hat{x_i},\ldots ,x_{j-1},[x_i,x_j]_{\mathfrak{g}},x_{j+1},\ldots, x_{n+1}),
 \end{align*}
where $f\in C^n_{LA}(\mathfrak{g},V)$ and $x_1,\ldots,x_{n+1} \in \mathfrak{g}$. 

 Then $\{C_{LA}^{n}(\mathfrak{g},V),\delta^n \}$ is a cochain complex. The corresponding cohomology groups are called the cohomology of $\mathfrak{g}$ with coefficients in the representation $V$ and the $n$th cohomology group is denoted by $H^n_{LA}(\mathfrak{g},V).$
 For details see Loday-Pirashvili cohomology  for Leibniz algebra in \cite{Loday}. We will  follow the notation $l_{V}(x,u)=[x,u] $ and $r_{V}(u,x)=[u,x]$ for all $x\in \mathfrak{g},~u \in V$  . Then the above coboundary map  $\delta ^n: C^n_{LA}(\mathfrak{g},V) \rightarrow C^{n+1}_{LA}(\mathfrak{g},V)$ becomes
 \begin{align*}
 (\delta^n(f))&(x_1,x_2,\ldots ,x_{n+1})\\
 ={ }&{ }\sum_{i=1}^{n}(-1)^{i+1}[x_i,f(x_1,\ldots,\hat{x_i},\ldots, x_{n+1})]+(-1)^{n+1}[f(x_1,\ldots,x_n),x_{n+1}] \\
 &+ \sum_{1 \leq i <j \leq n+1}(-1)^if(x_1,\ldots, \hat{x_i},\ldots ,x_{j-1},[x_i,x_j]_{\mathfrak{g}},x_{j+1},\ldots, x_{n+1}),
  \end{align*}
where $f\in C^n_{LA}(\mathfrak{g},V)$ and $x_1,\ldots,x_{n+1} \in \mathfrak{g}$.
Let $(\mathfrak{g}_T,[~,~]_{\mathfrak{g}})$ be a Rota-Baxter Leibniz algebra and $(V,l_V,r_V,T_V)$ be a representation of it. Now using proposition (3.1) and (3.2) we get a new Rota-Baxter Leibniz algebra $(\mathfrak{g}_T,[~,~]_*)$ with representation $(V,l^{'}_V,r_V^{'},T_V)$ induced by the Rota-Baxter operator. Now we consider the Loday-Pirashvili cochain complex of this induced Leibniz algebra $(\mathfrak{g},[~,~]_*)$ with representation $(V,l^{'}_V,r_V^{'})$
as follows:

For each  $n \geq 0$, we define cochain groups $C^n_{RBO}(\mathfrak{g},V):=\Hom(\mathfrak{g}^{\otimes n},V)$ and  boundary map $\partial ^n: C^n_{RBO}(\mathfrak{g},V) \rightarrow C^{n+1}_{RBO}(\mathfrak{g},V)$  by
\begin{align*}
(\partial^n&(f))(x_1,x_2,\ldots ,x_{n+1}) \\
={ }&{ }\sum_{i=1}^{n}(-1)^{i+1}l^{'}_{V}(x_i,f(x_1,\ldots,\hat{x_i},\ldots, x_{n+1}))+(-1)^{n+1}r^{'}_{V}(f(x_1,\ldots,x_n),x_{n+1})\\
& + \sum_{1 \leq i <j \leq n+1}(-1)^if(x_1,\ldots, \hat{x_i},\ldots ,x_{j-1},[x_i,x_j]_{*},x_{j+1},\ldots, x_{n+1}),\\
={ }&{ }\sum _{i=1}^{n}(-1)^{i+1}[T(x_i),f(x_1,\ldots,\hat{x_i},\ldots, x_{n+1})] -\sum _{i=1}^{n}(-1)^{i+1}T_V([x_i,f(x_1,\ldots,\hat{x_i},\ldots, x_{n+1})])\\
&+(-1)^{n+1}[f (x_1,\ldots ,x_{n}),T(x_{n+1})] -(-1)^{n+1}T_V([f (x_1,\ldots ,x_{n}),x_{n+1}])\\
&+\sum_{1\leq i< j\leq n+1}(-1)^if (x_1,\ldots,\hat{x_i},\ldots,x_{j-1},[T(x_i),x_j]_{\mathfrak{g}}+[x_i,T(x_j)]_{\mathfrak{g}},x_{j+1},\ldots ,x_{n+1})
 \end{align*}
where $f\in C^n_{RBO}(\mathfrak{g},V)$ and $x_1,\ldots,x_{n+1} \in \mathfrak{g}$.
Now one can observe that $\partial^{n+1} \circ \partial^n=0$. Hence, $\{C^{n}_{RBO}(\mathfrak{g},V),\partial^n\}$ is a cochain complex. This cochain complex is called the cochain complex of Rota-Baxter operator $T$ and  the corresponding cohomology groups are called the cohomology of Rota-Baxter operator $T$ with coefficients in the representation $V$ and is denoted by $H^n_{RBO}(\mathfrak{g},V).$

\begin{definition}
Let $(\mathfrak{g}_T,[~,~]_{\mathfrak{g}} )$ be a Rota-Baxter Leibniz algebra and $(V,l_V,r_V,T_V)$ be a representation of it. We define a map $\phi^{n} : C^{n}_{LA}(\mathfrak{g},V) \rightarrow C^{n}_{RBO}(\mathfrak{g},V)$ by 
\begin{align*}
\phi^n(f)&(x_1,x_2,\ldots,x_n)\\
={ }&{ }f(Tx_1,Tx_2,\ldots, Tx_n)- (T_v\circ f)(x_1,Tx_2,\ldots,Tx_n)\\
&-(T_v\circ f)(Tx_1,x_2,Tx_3,\ldots,Tx_n) -\ldots -(T_v\circ f)(Tx_1,Tx_2,\ldots,Tx_{n-1},x_n).
\end{align*}
\end{definition}
\begin{lemma}
For every $f\in C^n_{LA}(\mathfrak{g},V)$ and $x_1,\ldots,x_{n+1} \in \mathfrak{g}$, we have: 
\[
\phi^{n+1}(\delta^n(f))(x_1,x_2,x_3,\ldots,x_{n+1})
=\partial^n(\phi^n(f))(x_1,x_2,x_3,\ldots,x_{n+1}).
\]
\end{lemma}
\begin{proof}
Now,
\begin{align*}
&\phi^{n+1}(\delta^n(f))(x_1,x_2,x_3,\ldots,x_{n+1})\\
&=\delta^n(f)(Tx_1,Tx_2,Tx_3\ldots,Tx_{n+1})-(T_V\circ \delta^n(f))(x_1,Tx_2,Tx_3\ldots,Tx_{n+1})\\
&-(T_V\circ \delta^n(f))(Tx_1,x_2,Tx_3,\ldots,Tx_{n+1})-(T_V\circ \delta^n(f))(Tx_1,Tx_2,x_3,Tx_4,\ldots,Tx_{n+1})\\
&-\ldots-(T_V\circ \delta^n(f))(Tx_1,Tx_2,\ldots,Tx_n,x_{n+1})\\
&=\sum_{i=1}^n(-1)^{i+1}[Tx_i,f(Tx_1,Tx_2,\ldots,\widehat{Tx_i},\ldots,Tx_{n+1})]\\
&+(-1)^{n+1}[f(Tx_1,Tx_2,\ldots,Tx_n),Tx_{n+1}]\\
&+\sum_{1\leq i <j \leq n+1}(-1)^{i}f(Tx_1,\ldots,\widehat{Tx_i},\ldots,Tx_{j-1},[Tx_i,Tx_j]_{\mathfrak{g}},Tx_{j+1},\ldots,Tx_{n+1}) \\
&-T_V\bigg([x_1,f(Tx_2,Tx_3,\ldots , Tx_{n+1})]-[Tx_2,f(x_1,Tx_3,\ldots,Tx_{n+1})]\\
&+[Tx_3,f(x_1,Tx_2,Tx_4,\ldots,Tx_{n+1})]\\
&- \ldots +(-1)^{n+1}[Tx_n,f(x_1,Tx_2,\ldots,Tx_{n-1},Tx_{n+1})]\\
&+(-1)^{n+1}[f(x_1,Tx_2,Tx_3,\ldots,Tx_n),Tx_{n+1}]\\
& -\sum_{j=2}^{n+1} f(Tx_2,Tx_3,\ldots,Tx_{j-1},[x_1,Tx_j]_{\mathfrak{g}},Tx_j,\ldots,Tx_{n+1})\\
&+\sum_{2\leq i<j\leq n+1}(-1)^if(x_1,Tx_2,\ldots,\widehat{Tx_i},\ldots,Tx_{j-1},[Tx_i,Tx_j]_{\mathfrak{g}},Tx_{j+1},\ldots,Tx_{n+1}) \bigg)\\
&-T_V\bigg([Tx_1,f(x_2,Tx_3,\ldots , Tx_{n+1})]-[x_2,f(Tx_1,Tx_3,\ldots,Tx_{n+1})]\\
&+[Tx_3,f(Tx_1,x_2,Tx_4,\ldots,Tx_{n+1})]\\
&- \ldots +(-1)^{n+1}[Tx_n,f(Tx_1,x_2,Tx_3,\ldots,Tx_{n-1},Tx_{n+1})]\\
&+(-1)^{n+1}[f(Tx_1,x_2,Tx_3,\ldots,Tx_n),Tx_{n+1}]\\
& +\sum_{j=3}^{n+1} f(Tx_1,Tx_3,\ldots,Tx_{j-1},[x_2,Tx_j]_{\mathfrak{g}},Tx_{j+1},\ldots,Tx_{n+1})\\
&+\sum_{\substack{1\leq i<j\leq n+1 \\ i\neq 2}}(-1)^if(Tx_1,x_2,Tx_3,\ldots,\widehat{Tx_i},\ldots,Tx_{j-1},[Tx_i,Tx_j]_{\mathfrak{g}},Tx_{j+1},\ldots,Tx_{n+1}) \bigg)\\
&-T_V\bigg([Tx_1,f(Tx_2,x_3,Tx_4,\ldots , Tx_{n+1})]-[Tx_2,f(Tx_1,x_3,Tx_4,\ldots,Tx_{n+1})]\\
&+[x_3,f(Tx_1,Tx_2,Tx_4,\ldots,Tx_{n+1})]\\
&- \ldots +(-1)^{n+1}[Tx_n,f(Tx_1,Tx_2,x_3,Tx_4,\ldots,Tx_{n-1},Tx_{n+1})]\\
&+(-1)^{n+1}[f(Tx_1,Tx_2,x_3,Tx_4,\ldots,Tx_n),Tx_{n+1}]\\
& -\sum_{j=4}^{n+1} f(Tx_1,Tx_2,Tx_4,\ldots,Tx_{j-1},[x_3,Tx_j]_{\mathfrak{g}},Tx_{j+1},\ldots,Tx_{n+1})\\
&+\sum_{\substack{1\leq i<j\leq n+1 \\ i\neq 3}}(-1)^if(Tx_1,Tx_2,x_3,Tx_4\ldots,\widehat{Tx_i},\ldots,Tx_{j-1},[Tx_i,Tx_j]_{\mathfrak{g}},Tx_{j+1},\ldots,Tx_{n+1}) \bigg)\\
&-\ldots \\
&-T_V\bigg([Tx_1,f(Tx_2,Tx_3,\ldots ,Tx_n, x_{n+1})]-[Tx_2,f(Tx_1,Tx_3,\ldots,Tx_n,x_{n+1})]\\
&+[Tx_3,f(Tx_1,Tx_2,Tx_4,\ldots,Tx_n,x_{n+1})]- \ldots \\
&+(-1)^{n+1}[Tx_n,f(Tx_1,Tx_2,\ldots,Tx_{n-1},x_{n+1})]\\
&+(-1)^{n+1}[f(Tx_1,Tx_2,Tx_3,\ldots,Tx_n),x_{n+1}]\\
& +\sum_{i=1}^{n}(-1)^i f(Tx_1,Tx_2,\ldots,\widehat{Tx_i},\ldots,Tx_n,[Tx_i,x_{n+1}]_{\mathfrak{g}})\\
&+\sum_{\substack{1\leq i<j\leq n }}(-1)^if(Tx_1,Tx_2,\ldots,\widehat{Tx_i},\ldots,Tx_{j-1},[Tx_i,Tx_j]_{\mathfrak{g}},Tx_{j+1},\ldots,Tx_n,x_{n+1}) \bigg).
\end{align*}
Again,
\begin{align*}
&\partial^n(\phi^n(f))(x_1,x_2,x_3,\ldots,x_{n+1})\\
&=\sum _{i=1}^{n}(-1)^{i+1}[T(x_i),\phi^n(f)(x_1,\ldots,\hat{x_i},\ldots, x_{n+1})] \\
&-\sum _{i=1}^{n}(-1)^{i+1}T_V([x_i,\phi^n(f)(x_1,\ldots,\hat{x_i},\ldots, x_{n+1})])\\
&+(-1)^{n+1}[\phi^n(f) (x_1,\ldots ,x_{n}),T(x_{n+1})] -(-1)^{n+1}T_V([\phi^n(f) (x_1,\ldots ,x_{n}),x_{n+1}])\\
&+\sum_{1\leq i< j\leq n+1}(-1)^i \phi^n(f) (x_1,\ldots,\hat{x_i},\ldots,x_{j-1},[T(x_i),x_j]_{\mathfrak{g}}+[x_i,T(x_j)]_{\mathfrak{g}},x_{j+1},\ldots ,x_{n+1})\\
&=\bigg([Tx_1,\phi^n(f)(x_2,x_3,\ldots,x_{n+1})]-[Tx_2,\phi^n(f)(x_1,x_3,\ldots,x_{n+1})]\\
&+[Tx_3,\phi^n(f)(x_1,x_2,x_4,\ldots,x_{n+1})]\\
&-\ldots +(-1)^{n+1}[Tx_n,\phi^n(f)(x_1,x_2,\ldots,x_{n-1},x_{n+1})]\bigg)\\
&-\bigg(T_V([x_1,\phi^n(f)(x_2,x_3,\ldots,x_{n+1})])-T_V([x_2,\phi^n(f)(x_1,x_3,\ldots,x_{n+1})])\\
&+T_V([x_3,\phi^n(f)(x_1,x_2,x_4,\ldots,x_{n+1})])\\
&-\ldots +(-1)^{n+1}T_V([x_n,\phi^n(f)(x_1,x_2,\ldots,x_{n-1},x_{n+1})])\bigg)\\
&+(-1)^{n+1}[\phi^n(f) (x_1,\ldots ,x_{n}),T(x_{n+1})] -(-1)^{n+1}T_V([\phi^n(f) (x_1,\ldots ,x_{n}),x_{n+1}])\\
&-\sum_{j=2}^{n+1}\phi^n(f)(x_2,x_3,\ldots,x_{j-1},[Tx_1,x_j]_{\mathfrak{g}}+[x_1,Tx_j]_{\mathfrak{g}},x_{j+1},\ldots,x_{n+1})\\
&+\sum_{j=3}^{n+1}\phi^n(f)(x_1,x_3,\ldots,x_{j-1},[Tx_2,x_j]_{\mathfrak{g}}+[x_2,Tx_j]_{\mathfrak{g}},x_{j+1},\ldots,x_{n+1}) \\
&-\sum_{j=4}^{n+1}\phi^n(f)(x_1,x_2,x_4,\ldots,x_{j-1},[Tx_3,x_j]_{\mathfrak{g}}+[x_3,Tx_j]_{\mathfrak{g}},x_{j+1},\ldots,x_{n+1})\\
&+\ldots \\
&+(-1)^n\sum_{j=n+1}^{n+1}\phi^n(f)(x_1,x_2,x_3,\ldots,x_{j-1},[Tx_n,x_j]_{\mathfrak{g}}+[x_n,Tx_j]_{\mathfrak{g}},x_{j+1},\ldots,x_{n-1},x_{n+1})\\
&=\bigg([Tx_1,f(Tx_2,Tx_3,\ldots,Tx_{n+1})]-[Tx_1,T_V(f(x_2,Tx_3,\ldots,Tx_{n+1}))]\\
&-[Tx_1,T_V(f(Tx_2,x_3,Tx_4,\ldots,Tx_{n+1}))]\\
&-[Tx_1,T_V(f(Tx_2,Tx_3,x_4,Tx_5,\ldots,Tx_{n+1}))]\\
&- \ldots-[Tx_1,T_V(f(Tx_2,Tx_3,\ldots,Tx_{n},x_{n+1}))]\bigg)\\
&-\bigg([Tx_2,f(Tx_1,Tx_3,\ldots,Tx_{n+1})]-[Tx_2,T_V(f(x_1,Tx_3,\ldots,Tx_{n+1}))]\\
&-[Tx_2,T_V(f(Tx_1,x_3,Tx_4,\ldots,Tx_{n+1}))]\\
&-[Tx_2,T_V(f(Tx_1,Tx_3,x_4,Tx_5,\ldots,Tx_{n+1}))]\\
&- \ldots-[Tx_2,T_V(f(Tx_1,Tx_3,\ldots,Tx_{n},x_{n+1}))]\bigg)\\
&+\bigg([Tx_3,f(Tx_1,Tx_2,Tx_4,\ldots,Tx_{n+1})]-[Tx_3,T_V(f(x_1,Tx_2,Tx_4,\ldots,Tx_{n+1}))]\\
&-[Tx_3,T_V(f(Tx_1,x_2,Tx_4\ldots,Tx_{n+1}))]\\
&-[Tx_3,T_V(f(Tx_1,Tx_2,x_4,Tx_5,\ldots,Tx_{n+1}))]\\
&- \ldots-[Tx_3,T_V(f(Tx_1,Tx_2,Tx_4,\ldots,Tx_{n},x_{n+1}))]\bigg)\\
&-\ldots +\\
&(-1)^{n+1}\bigg([Tx_n,f(Tx_1,Tx_2,\ldots,Tx_{n-1},Tx_{n+1})]-[Tx_n,T_V(f(x_1,Tx_2,\ldots,Tx_{n-1},Tx_{n+1}))]\\
&-[Tx_n,T_V(f(Tx_1,x_2,Tx_3,\ldots,Tx_{n-1},Tx_{n+1}))]\\
&-[Tx_n,T_V(f(Tx_1,Tx_2,x_3,Tx_4,\ldots,Tx_{n-1},Tx_{n+1}))]- \ldots \\
&-[Tx_n,T_V(f(Tx_1,Tx_2,Tx_3,\ldots,Tx_{n-1},x_{n+1}))]\bigg)\\
&-T_V\bigg([x_1,f(Tx_2,Tx_3,\ldots,Tx_{n+1})]-[x_1,T_V(f(x_2,Tx_3,\ldots,Tx_{n+1}))]\\
&-[x_1,T_V(f(Tx_2,x_3,Tx_4,\ldots,Tx_{n+1}))] \\
&-[x_1,T_V(f(Tx_2,Tx_3,x_4,Tx_5,\ldots,Tx_{n+1}))]- \ldots-[x_1,T_V(f(Tx_2,Tx_3,\ldots,Tx_{n},x_{n+1}))]\bigg)\\
&+T_V\bigg([x_2,f(Tx_1,Tx_3,\ldots,Tx_{n+1})]-[x_2,T_V(f(x_1,Tx_3,\ldots,Tx_{n+1}))]\\
&-[x_2,T_V(f(Tx_1,x_3,Tx_4,\ldots,Tx_{n+1}))]\\
&-[x_2,T_V(f(Tx_1,Tx_3,x_4,Tx_5,\ldots,Tx_{n+1}))]- \ldots-[x_2,T_V(f(Tx_1,Tx_3,\ldots,Tx_{n},x_{n+1}))]\bigg)\\
&-T_V\bigg([x_3,f(Tx_1,Tx_2,Tx_4,\ldots,Tx_{n+1})]-[x_3,T_V(f(x_1,Tx_2,Tx_4,\ldots,Tx_{n+1}))]\\
&-[x_3,T_V(f(Tx_1,x_2,Tx_4\ldots,Tx_{n+1}))]\\
&-[x_3,T_V(f(Tx_1,Tx_2,x_4,Tx_5,\ldots,Tx_{n+1}))]\\
&- \ldots-[x_3,T_V(f(Tx_1,Tx_2,Tx_4,\ldots,Tx_{n},x_{n+1}))]\bigg)\\
&+\ldots \\ 
&-(-1)^{n+1}T_V\bigg([x_n,f(Tx_1,Tx_2,\ldots,Tx_{n-1},Tx_{n+1})]\\
&-[x_n,T_V(f(x_1,Tx_2,\ldots,Tx_{n-1},Tx_{n+1}))]\\
&-[x_n,T_V(f(Tx_1,x_2,Tx_3,\ldots,Tx_{n-1},Tx_{n+1}))]\\
&-[x_n,T_V(f(Tx_1,Tx_2,x_3,Tx_4,\ldots,Tx_{n-1},Tx_{n+1}))]\\
&- \ldots-[x_n,T_V(f(Tx_1,Tx_2,Tx_3,\ldots,Tx_{n-1},x_{n+1}))]\bigg)\\
&+(-1)^{n+1}\bigg([f(Tx_1,Tx_2,Tx_3,\ldots,Tx_n),Tx_{n+1}]-[f(x_1,Tx_2,Tx_3,\ldots,Tx_n),Tx_{n+1}]\\
&-[f(Tx_1,x_2,Tx_3,\ldots,Tx_n),Tx_{n+1}]-[f(Tx_1,Tx_2,x_3,Tx_4,\ldots,Tx_n),Tx_{n+1}]\\
&-\ldots-[f(Tx_1,Tx_2,Tx_3,\ldots,Tx_{n-1},x_n),Tx_{n+1}]\bigg)\\
&-(-1)^{n+1}T_V\bigg([f(Tx_1,Tx_2,Tx_3,\ldots,Tx_n),x_{n+1}]-[f(x_1,Tx_2,Tx_3,\ldots,Tx_n),x_{n+1}]\\
&-[f(Tx_1,x_2,Tx_3,\ldots,Tx_n),Tx_{n+1}]-[f(Tx_1,Tx_2,x_3,Tx_4,\ldots,Tx_n),x_{n+1}]\\
&-\ldots-[f(Tx_1,Tx_2,Tx_3,\ldots,Tx_{n-1},x_n),x_{n+1}]\bigg)\\
&-\sum_{j=2}^{n+1}\bigg(f(Tx_2,Tx_3,\ldots,Tx_{j-1},[Tx_1,Tx_j]_{\mathfrak{g}},Tx_{j+1},\ldots,Tx_{n+1})\\
&-T_V(f(x_2,Tx_3,\ldots,Tx_{j-1},[Tx_1,Tx_j]_{\mathfrak{g}},Tx_{j+1},\ldots,Tx_{n+1})) \\
&-T_V(f(Tx_2,x_3,Tx_4,\ldots,Tx_{j-1},[Tx_1,Tx_j]_{\mathfrak{g}},Tx_{j+1},\ldots,Tx_{n+1}))\\
&-\ldots -T_V(f(Tx_2,Tx_3,\ldots,Tx_{j-2},x_{j-1},[Tx_1,Tx_j]_{\mathfrak{g}},Tx_{j+1},\ldots,Tx_{n+1}))\\
&-T_V(f(Tx_2,Tx_3,\ldots,Tx_{j-1},[x_1,Tx_j]_{\mathfrak{g}}+[Tx_1,x_j]_{\mathfrak{g}},Tx_{j+1},\ldots,Tx_{n+1}))-\ldots\\
&-T_V(f(Tx_2,Tx_3,\ldots,Tx_{j-1},[Tx_1,Tx_j]_{\mathfrak{g}},Tx_{j+1},\ldots,Tx_{n},x_{n+1}))\bigg)\\
&+\sum_{j=3}^{n+1}\bigg(f(Tx_1,Tx_3,\ldots,Tx_{j-1},[Tx_2,Tx_j]_{\mathfrak{g}},Tx_{j+1},\ldots,Tx_{n+1})\\
&-T_V(f(x_1,Tx_3,\ldots,Tx_{j-1},[Tx_2,Tx_j]_{\mathfrak{g}},Tx_{j+1},\ldots,Tx_{n+1}))\\
&-T_V(f(Tx_1,x_3,Tx_4\ldots,Tx_{j-1},[Tx_2,Tx_j]_{\mathfrak{g}},Tx_{j+1},\ldots,Tx_{n+1}))\\
&-\ldots - T_V(f(Tx_1,Tx_3,\ldots,Tx_{j-2},x_{j-1},[Tx_2,Tx_j]_{\mathfrak{g}},Tx_j,\ldots,Tx_{n+1}))\\
&-T_V(f(Tx_1,Tx_3,\ldots,Tx_{j-1},[x_2,Tx_j]_{\mathfrak{g}}+[Tx_2,x_j]_{\mathfrak{g}},Tx_{j+1},\ldots,Tx_{n+1}))\\
&- \ldots -T_V(f(Tx_1,Tx_3,\ldots,Tx_{j-1},[Tx_2,Tx_j]_{\mathfrak{g}},Tx_j,\ldots,Tx_{n},x_{n+1}))\bigg)\\
&-\sum_{j=4}^{n+1}\bigg(f(Tx_1,Tx_2,Tx_4,\ldots,Tx_{j-1},[Tx_3,Tx_j]_{\mathfrak{g}},Tx_{j+1},\ldots,Tx_{n+1})\\
&-T_V(f(x_1,Tx_2,Tx_4,\ldots,Tx_{j-1},[Tx_3,Tx_j]_{\mathfrak{g}},Tx_{j+1},\ldots,Tx_{n+1}))\\
&-T_V(f(Tx_1,x_2,Tx_4,\ldots,Tx_{j-1},[Tx_3,Tx_j]_{\mathfrak{g}},Tx_{j+1},\ldots,Tx_{n+1}))\\
&- T_V(f(Tx_1,Tx_2,x_4,Tx_5\ldots,x_{j-1},[Tx_3,Tx_j]_{\mathfrak{g}},Tx_{j+1},\ldots,Tx_{n+1}))\\
&- \ldots -T_V(f(Tx_1,Tx_2,Tx_4,\ldots,Tx_{j-1},[x_3,Tx_j]_{\mathfrak{g}}+[Tx_3,x_j]_{\mathfrak{g}},Tx_{j+1},\ldots,Tx_{n+1}))-\ldots\\
&-T_V(f(Tx_1,Tx_2,Tx_4,\ldots,Tx_{j-1},[Tx_3,Tx_j]_{\mathfrak{g}},Tx_{j+1},\ldots,Tx_{n},x_{n+1}))\bigg)\\
&+\ldots \\
&+(-1)^n \sum_{j=n+1}^{n+1}\bigg(f(Tx_1,Tx_2,\ldots,Tx_{j-1},[Tx_n,Tx_j]_{\mathfrak{g}},Tx_{j+1},\ldots,Tx_{n-1},Tx_{n+1})\\
&-T_V(f(x_1,Tx_2,Tx_3,\ldots,Tx_{j-1},[Tx_n,Tx_j]_{\mathfrak{g}},Tx_{j+1},\ldots,Tx_{n-1},Tx_{n+1}))\\
&-T_V(f(Tx_1,x_2,Tx_3,\ldots,Tx_{j-1},[Tx_n,Tx_j]_{\mathfrak{g}},Tx_{j+1},\ldots,Tx_{n-1},Tx_{n+1}))\\
&-\ldots T_V(f(Tx_1,Tx_2,x_3,Tx_4\ldots,Tx_{j-1},[Tx_n,Tx_j]_{\mathfrak{g}},Tx_{j+1},\ldots,Tx_{n-1},Tx_{n+1}))\\
&- \ldots -T_V(f(Tx_1,Tx_2,Tx_3,\ldots,Tx_{j-1},[x_n,Tx_j]_{\mathfrak{g}}+[Tx_n,x_j]_{\mathfrak{g}},Tx_{j+1},\ldots,Tx_{n-1},Tx_{n+1}))\\
&-\ldots-T_V(f(Tx_1,Tx_2,Tx_3,\ldots,Tx_{j-1},[Tx_n,Tx_j]_{\mathfrak{g}},Tx_{j+1},\ldots,Tx_{n-1},x_{n+1}))\bigg).
\end{align*}
Now, using equations 
\[[T(x),T_V(u)]=T_V(([T(x),u]+[x,T_V(u)])\]
\[[T_V(u),T(x)]=T_V([T_V(u),x]+[u,T(x)])\]
for all $x \in \mathfrak{g}$ and $u \in V$ 
, we get
\[
\phi^{n+1}(\delta^n(f))(x_1,x_2,x_3,\ldots,x_{n+1})=\partial^n(\phi^n(f))(x_1,x_2,x_3,\ldots,x_{n+1}).
\qedhere
\]
\end{proof}
Now by the lemma (4.2), we have the following commutative diagram.
\[\begin{tikzcd}
	{C^1_{LA}(\mathfrak{g},V)} & {C^2_{LA}(\mathfrak{g},V)} & {C^n_{LA}(\mathfrak{g},V)} & {C^{n+1}_{LA}(\mathfrak{g},V)} & {} \\
	{C^1_{RBO}(\mathfrak{g},V)} & {C^2_{RBO}(\mathfrak{g},V)} & {C^n_{RBO}(\mathfrak{g},V)} & {C^{n+1}_{RBO}(\mathfrak{g},V)} & {}
	\arrow[from=1-1, to=1-2,]{r}{\delta^1}
	\arrow[from=1-1, to=2-1]{d}{\phi^1}
	\arrow[from=1-2, to=2-2]{d}{\phi^2}
	\arrow[from=2-1, to=2-2]{r}{\partial^1}
	\arrow[dotted, no head, from=1-2, to=1-3]
	\arrow[dotted, no head, from=2-2, to=2-3]
	\arrow[from=1-3, to=1-4]{r}{\delta^n}
	\arrow[from=1-3, to=2-3]{d}{\phi^n}
	\arrow[from=2-3, to=2-4]{d}{\partial^n}
	\arrow[from=1-4, to=2-4]{r}{\phi^{n+1}}
	\arrow[dotted, no head, from=1-4, to=1-5]
	\arrow[dotted, no head, from=2-4, to=2-5].
\end{tikzcd}\]

Now, we combine the cochain complex of Leibniz algebra and the cochain complex of Rota-Baxter operator to define the cochain complex of Rota-Baxter Leibniz algebra.

Let $(\mathfrak{g}_T,[~,~]_{\mathfrak{g}} )$ be a Rota-Baxter Leibniz algebra and $(V,l_V,r_V,T_V)$ be a representation of it. Now we define the cochain groups by
\[C^{0}_{RBLA}(\mathfrak{g},V)= C^{0}_{LA}(\mathfrak{g},V) ~~ \mbox{and}~~ C^n_{RBLA}(\mathfrak{g},V)= C^n_{LA}(\mathfrak{g},V)\oplus C_{RBO}^{n-1}(\mathfrak{g},V), \forall n \geq 1,\] and the coboundary map
$d^n: C^n_{RBLA}(\mathfrak{g},V) \rightarrow C^{n+1}_{RBLA}(\mathfrak{g},V)$ is defined by
\[d^n (\alpha , \beta)=(\delta^n(\alpha),- \partial^{n-1}(\beta)-\phi^n(\alpha) )\]
for any $\alpha \in C^n_{LA}(\mathfrak{g},V)$ and $\beta \in C^{n-1}_{RBO}(\mathfrak{g},V).$

\begin{theorem}
The map $d^n: C^n_{RBLA}(\mathfrak{g},V) \rightarrow C^{n+1}_{RBLA}(\mathfrak{g},V)$ satisfies $d^{n+1}\circ d^n=0$.
\end{theorem}
\begin{proof}
Let $f\in C^n_{LA}(\mathfrak{g},V) $ and $g\in C^{n-1}_{RBO}(\mathfrak{g},V)$, then we have  
\begin{align*}
d^{n+1}\circ d^{n}(f,g)&= d^{n+1}(\delta ^n(f),-\partial ^{n-1}(g)-\phi ^{n}(f))\\
&=( \delta ^{n+1}(\delta ^n(f)),-\partial^{n}(-\partial^{n-1}(g)
-\phi^n(f))-\phi^{n+1}(\delta ^n(f)))\\
&=(0,~\partial^n(\phi^n(f))-\phi^{n+1}(\delta^n(f)))=0.
\qedhere
\end{align*}
\end{proof}
Therefore, it follows from the above theorem that $\{C^n_{RBLA}(\mathfrak{g},V),d^n\}$ is a cochain complex and the corresponding cohomology groups are called \textbf{cohomology of Rota-Baxter Leibniz algebra} $(\mathfrak{g}_T,[~,~]_{\mathfrak{g}} )$ with coefficients in the representation $V$  and is denoted by $\mathbf{H^n_{RBLA}(\mathfrak{g},V)}$, $~n\geq 0.$ Note that  in this case there exists a short exact sequence of complexes 
\[0 \longrightarrow C^n_{RBO}(\mathfrak{g},V) \longrightarrow C^n_{RBLA}(\mathfrak{g},V) \longrightarrow C^n_{LA}(\mathfrak{g},V) \longrightarrow 0.\]

\section{Deformations of Rota-Baxter Leibniz algebras}\label{sec4}
In this section, we study a one-parameter formal deformation of Rota-Baxter Leibniz algebra. We denote the bracket $[~,~]_{\mathfrak{g}}$ by $\mu.$

\begin{definition}
A formal one-parameter deformation of a Rota-Baxter Leibniz algebra $(\mathfrak{g}_T,\mu)$ is a pair of two power series  $(\mu_t,T_t)$
\[ \mu_t=\sum _{i=0}^{\infty}\mu_it^i, ~ \mu_{i} \in C^2_{LA}(\mathfrak{g},\mathfrak{g}),~~~~ T_t= \sum_{i=0}^{\infty} T_it^i,~ T_i \in C^1_{RBO}(\mathfrak{g},\mathfrak{g}),\]
 such that $(\mathfrak{g}[[t]]_{T_t},\mu_t)$ is a Rota-Baxter Leibniz algebra with $(\mu_0,T_{0})=(\mu ,T)$, where $\mathfrak{g}[[t]]$, the space of formal power series in $t$ with coefficients from $\mathfrak{g}$ is a $\mathbb{K}[[t]]$ module, $\mathbb{K}$ being the ground field of $(\mathfrak{g}_T,\mu)$.
\end{definition}
The above definition holds if and only if for any $x,y,z \in \mathfrak{g}$ the following conditions are satisfied 
\[\mu_t(x,\mu_t(y,z))=\mu_t(\mu_t(x,y),z)+\mu_t(y,\mu_t(x,z)),\]
and \[ \mu_t(T_t(x),T_t(y))=T_t(\mu_t(x,T_t(y))+\mu_t(T_t(x),y)).\]

Expanding the above equations and equating the coefficients of $t^n$ from both sides we have 
\begin{align}
\sum _{\substack{i+j=n \\i,j\geq 0}} \mu_i(x,\mu_j (y,z))=\sum _{\substack{i+j=n \\i,j\geq 0}}\mu_i(\mu _j(x,y),z)+\sum _{\substack{i+j=n \\i,j\geq 0}} \mu_i(y, \mu_j (x,z)),
\end{align}
and 
\begin{align}
\sum_{\substack{i+j+k=n \\ i,j,k \geq 0}} \mu_i(T_j(u),T_k(v))=\sum_{\substack{i+j+k=n \\ i,j,k \geq 0}}T_i(\mu_j(T_k(u),v))+\sum_{\substack{i+j+k=n \\ i,j,k \geq 0}}T_i(\mu_j(u,T_k(v))).
\end{align}
Observe that for $n=0$, the above conditions are exactly the conditions in the definitions of  Leibniz algebra and the Rota-Baxter operator.
\begin{definition}
The infinitesimal of the deformation $(\mu_t, T_t)$ is the pair $(\mu_1, T_1)$. Suppose more  generally that $(\mu_n, T_n)$ is the first non-zero term of $(\mu_t, T_t)$ after $(\mu_0, T_0)$, such $(\mu_n, T_n)$ is called a $n$-infinitesimal of the deformation.
\end{definition}
\begin{theorem}
Let $(\mu_t,T_t)$ be a formal one-parameter deformation of Rota-Baxter Leibniz algebra $(\mathfrak{g}_T,\mu)$. Then $(\mu_1,T_1)$ is a $2$-cocycle in the cochain complex $\{C^{n}_{RBLA}(\mathfrak{g},\mathfrak{g}),d^n\}.$
\end{theorem}
\begin{proof}
Putting $n=1$ in the equation (5.1) we get 
\[ \mu (x,\mu_1(y,z))+\mu_1(x,\mu (y,z))= \mu (\mu_1(x,y),z)+\mu_1(\mu (x,y),z) \\
+\mu_1 (y, \mu (x,z))+\mu (y, \mu_1(x,z)).\]
This gives
$\delta^2(\mu_1)(x,y,z)=0 \in C^{2}_{LA}(\mathfrak{g},\mathfrak{g})$.
Again, putting $n=1$ in (5.2) we get 
\begin{align*}
\mu_1&(T(x_1),T(x_2))+\mu (T_1(x_1),T(x_2))+\mu (T(x_1),T_1(x_2)) \\
&-T_1(\mu (T(x_1),x_2)) -T(\mu (T_1(x_1),x_2)) -T(\mu _1 (T(x_1),x_2)) \\
&-T_1(\mu (x_1,T(x_2)))-T(\mu _1 (x_1,T(x_2)))-T(\mu (x_1,T_1(x_2))) \
= \ 0.
\end{align*}
This gives
\begin{align*}
-\partial ^1(T_1)(x_1,x_2)&=  -T(\mu_1 (x_1,T(x_2)))-T(\mu_1(T(x_1),x_2))+\mu_1(T(x_1),T(x_2))\\
                  &= \phi ^2(\mu_1)(x_1,x_2).
\end{align*}
Therefore, $-\partial^(T_1)-\phi^2(\mu_1)=0$.
Hence, $d^2(\mu_1,T_1)=0$. Thus, $(\mu_1,T_1)$ is a $2$-cocycle in the cochain complex $\{C^{n}_{RBLA}(\mathfrak{g},\mathfrak{g}),d^n\}.$
\end{proof}

\begin{theorem}
Let $(\mu_t,T_t)$ be a formal one-parameter deformation of Rota-Baxter Leibniz algebra $(\mathfrak{g}_T,\mu)$. Then $n$-infinitesimal of the deformation is a $2$-cocycle.
\end{theorem}
\begin{proof}
The proof is similar to the above theorem.
\end{proof}

\begin{definition}
Let $(\mu_t,T_t)$ and $(\mu_t^{'},T_t^{'})$ be two formal one-parameter deformations of a Rota-Baxter Leibniz algebra $(\mathfrak{g}_T,\mu )$. A formal isomorphism from 
$(\mu_t,T_t)$ to $(\mu_t^{'},T_t^{'})$ is a power series $\psi _t=\sum _{i=0}\psi_i t^i : \mathfrak{g}[[t]] \rightarrow \mathfrak{g}[[t]]$, where $\psi_i: \mathfrak{g} \rightarrow \mathfrak{g}$ are linear maps with $\psi_0$ is the identity map on $\mathfrak{g}$ and also the following conditions are satisfied.
\begin{align}
&\psi_t \circ \mu^{'}_t=\mu_t \circ (\psi_t \otimes \psi_t)\\
& \psi_t \circ T_t^{'}=T_{t} \circ \psi_t.  
\end{align}
In this case, we say that $(\mu_t,T_t)$ and $(\mu_t^{'},T_t^{'})$  are equivalent.
Note that the equation (5.3) and (5.4) can be written as follows respectively:
\begin{align}
&\sum_{\substack {i+j=n \\ i,j\geq 0}}\psi _i(\mu_j^{'}(x,y))=\sum_{\substack {i+j+k=n \\ i,j,k\geq 0}}\mu_i(\psi_j(x),\psi_k(y)),~~ x,y \in \mathfrak{g},\\
& \sum_{\substack {i+j=n \\ i,j\geq 0}}\psi _i \circ T^{'}_j=\sum_{\substack {i+j=n \\ i,j\geq 0}} T_i \circ \psi _j.
\end{align}
\end{definition}
\begin{theorem}
The infinitesimal of two equivalent formal one-parameter deformations of Rota-Baxter Leibniz algebra $(\mathfrak{g}_T , \mu )$ is in the same cohomology class. 
\end{theorem}
\begin{proof}
Let $\psi_t : (\mu_t,T_t) \rightarrow  (\mu_t^{'},T_t^{'})$ be a formal isomorphism. Now putting $n=1$ in equation (5.5) and (5.6) we get 
\begin{align*}
&\mu^{'}_1(x,y)=\mu_1(x,y)+\mu (x,\psi_1 (y))+\mu (\psi_1(x),y)-\psi_1(\mu (x,y)) ,~~ x,y \in \mathfrak{g}\\
& T_1^{'}=T_1+T \circ \psi_1-\psi _1 \circ T
\end{align*}
Therefore, we have
\[(\mu_1^{'},T_1^{'})-(\mu_1,T_1)=(\delta^1(\psi_1),-\phi^1(\psi_1))=d^1(\psi_1,0) \in C^{1}_{RBLA}(\mathfrak{g},\mathfrak{g}).\qedhere\]
\end{proof}
\begin{definition}
A Rota-Baxter Leibniz algebra is called rigid if every formal one-parameter deformation is trivial.
\end{definition}
\begin{theorem}
Let $(\mathfrak{g}_T,\mu)$ be a Rota-Baxter Leibniz algebra. If $H^2_{RBLA}(\mathfrak{g},\mathfrak{g})=0$, then $(\mathfrak{g}_T,\mu)$ is rigid.
\end{theorem}
\begin{proof}
Let $(\mu_t,T_t)$ be a formal one-parameter deformation of $(\mathfrak{g}_T,\mu)$. Since $(\mu_1,T_1)$ is a $2$-cocycle and $H^2_{RBLA}(\mathfrak{g},\mathfrak{g})=0$, thus, there exists a map $\psi_1^{'}$ and $x\in \mathbb{K} $, where $\mathbb{K}$ is the ground field of  Rota-Baxter Leibniz algebra $(\mathfrak{g}_T,\mu)$, such that
\[(\psi_1^{'},x) \in C^1_{RBLA}(\mathfrak{g},\mathfrak{g})=C^1_{LA}(\mathfrak{g},\mathfrak{g})\oplus \Hom(\mathbb{K},\mathfrak{g})\]
and $(\mu_1,T_1)=d^1(\psi_1^{'},x)$. Hence, $\mu_1=\delta^1(\psi_1^{'})~ \mbox{and}~T_1=-\partial^0(x)-\phi^1(\psi_1)$. If $\psi_1=\psi_1^{'}+\delta^0(x)$, then $\mu_1=\delta^1(\psi_1),~ T_1=-\phi^1(\psi_1)$. Now, let $\psi_t=Id_{\mathfrak{g}}-t\psi_t$. Then we have two equivalent deformation 
$(\mu_t,T_t)$ and $(\bar{\mu_t},\bar{T_t})$, where
\[\bar{\mu_t}=\psi_t^{-1} \circ \mu_t \circ (\psi_t \times \psi_t), ~~ \bar{T_t}=\psi_t^{-1} \circ T_t \circ \psi_t.
\]
Now by theorem (5.4) we have,  $\bar{\mu_1}=0,\bar{T_1}=0$. Hence,
\begin{align*}
&\bar{\mu_t}=\mu +\bar{\mu_2}t^2+\ldots ,\\
& \bar{T_t}=T+\bar{T_2}t^2+\ldots
\end{align*}
Thus, the linear terms of  $(\bar{\mu_2},\bar{T_2})$ vanishes, hence, repeatedly applying the same argument we conclude that $(\mu_t,T_t)$ is equivalent to the trivial  deformation. Hence, $(\mathfrak{g}_T,\mu)$ is rigid. 
\end{proof}

\section{Abelian extensions of Rota-Baxter Leibniz algebras}\label{sec5}
Let $(\mathfrak{g}_{T}, [~,~]_{\mathfrak{g}})$ be a Rota-Baxter Leibniz algebra and $V$ be a vector space. Observe that if $T_V$ is a linear operator on the vector space $V$ and  if we define the bracket by $\mu(x,y)=0$ for all $x,y \in V$. Then $(V_{T_V},\mu)$ has a structure of Rota-Baxter Leibniz algebra.
\begin{definition}
 An abelian extension of the Rota-Baxter Leibniz algebra $(\mathfrak{g}_{T},[~,~]_{\mathfrak{g}})$ is a short exact sequence of morphisms of Rota-Baxter Leibniz algebra 
 \[
\begin{tikzcd}
0 \arrow[r] & (V_{T_V},\mu) \arrow[r ,"i"] & (\hat{\mathfrak{g}}_{\hat{T}},[~,~]_{\wedge}) \arrow[r,"p"] & (\mathfrak{g}_{T},[~,~]_{\mathfrak{g}}) \arrow [r] & 0 
\end{tikzcd} ,
\]
that is, there exists a commutative diagram 
\[
\begin{tikzcd}
0 \arrow[r] & V \arrow[r ,"i"] \arrow[d,"T_V"]& \hat{\mathfrak{g}} \arrow[d,"\hat{T}"]\arrow[r,"p"] & \mathfrak{g} \arrow [r]\arrow[d,"T"]  & 0 \\
 0 \arrow[r] & V \arrow[r ,"i"] & \hat{\mathfrak{g}} \arrow [r,"p"]  &  \mathfrak{g} \arrow [r]  & 0
\end{tikzcd}
\]
where $\mu (a,b)=0$ for all $a,b \in V.$ In this case we say that $(\hat{\mathfrak{g}}_{\hat{T}},[~,~]_{\wedge})$ is an abelian extension of the Rota-Baxter Leibniz algebra $(\mathfrak{g}_{T},[~,~]_{\mathfrak{g}})$ by $(V_{T_V},\mu).$
 \end{definition}

\begin{definition}
Let $(\hat{\mathfrak{g}}_{\hat{T_1}},[~,~]_{\wedge_1})$ and $(\hat{\mathfrak{g}}_{\hat{T_2}},[~,~]_{\wedge_2})$ be two abelian extension of $(\mathfrak{g}_{T},[~,~]_{\mathfrak{g}})$ by $(V_{T_V},\mu)$. Then this two extension are said to be isomorphic if there exists an isomorphism of Rota-Baxter Leibniz algebra $\xi : (\hat{\mathfrak{g}}_{\hat{T_1}},[~,~]_{\wedge_1}) \rightarrow  (\hat{\mathfrak{g}}_{\hat{T_2}},[~,~]_{\wedge_2}) $ so that the following diagram is commutative : 
\[
\begin{tikzcd}
0 \arrow[r] & (V_{T_V},\mu) \arrow[r ,"i"] \arrow[d,equal]& (\hat{\mathfrak{g}}_{\hat{T_1}},[~,~]_{\wedge_1})\arrow[d,"\xi"]\arrow[r,"p"] & (\mathfrak{g}_T,[~,~]_{\mathfrak{g}}) \arrow [r]\arrow[d,equal]  & 0 \\
 0 \arrow[r] & (V_{T_V},\mu) \arrow[r ,"i"] & (\hat{\mathfrak{g}}_{\hat{T_2}},[~,~]_{\wedge_2}) \arrow [r,"p"]  &  (\mathfrak{g}_T,[~,~]_{\mathfrak{g}}) \arrow [r]  & 0.
\end{tikzcd}
\]
\end{definition}

\begin{definition}
A section of an abelian extension $(\hat{\mathfrak{g}}_{\hat{T}},[~,~]_{\wedge})$ of $(\mathfrak{g}_T,[~,~]_{\mathfrak{g}})$ by $(V_{T_V}, \mu)$ is a linear map $s : \mathfrak{g} \rightarrow \hat{\mathfrak{g}}$ such that $p \circ s= Id_{\mathfrak{g}}.$
 \end{definition}
\begin{definition}
Let $(\hat{\mathfrak{g}}_{\hat{T}},[~,~]_{\wedge})$ be an abelian extension of $(\mathfrak{g}_T,[~,~]_{\mathfrak{g}})$ by $(V_{T_V},\mu)$ with a section $s: \mathfrak{g} \rightarrow \hat{\mathfrak{g}}$. Now define $\bar{l}_V: \mathfrak{g} \otimes V \rightarrow V$ and $\bar{r}_V: V \otimes \mathfrak{g} \rightarrow V $ by respectively  $\bar{l}_V(x,u)=[s(x),u]_{\wedge}$ and $\bar{r}_V(u,x)=[u,s(x)]_{\wedge}$ for all $x\in \mathfrak{g}, u \in V.$
\end{definition}
\begin{theorem}
Let $(\hat{\mathfrak{g}}_{\hat{T}},[~,~]_{\wedge})$ be an abelian extension of $(\mathfrak{g}_T,[~,~]_{\mathfrak{g}})$ by $(V_{T_V},\mu)$ with a section $s: \mathfrak{g} \rightarrow \hat{\mathfrak{g}}$. Then, $(V,\bar{l}_V,\bar{r}_V,T_V)$ is a representation  Rota-Baxter Leibniz algebra  $(\mathfrak{g}_T,[~,~]_{\mathfrak{g}}).$
\end{theorem}
\begin{proof}
As $s([x,y]_{\mathfrak{g}})-[s(x),s(y)]_{\wedge} \in V$, hence, $[s([x,y]_{\mathfrak{g}}),u]_{\wedge}=[[s(x),s(y)]_{\wedge},u]_{\wedge}$ for all $x,y \in \mathfrak{g}, u \in V$. Therefore,  we have
\begin{align*}
\bar{l}_V&(x,\bar{l}_V(y,u))-\bar{l}_V([x,y]_{\mathfrak{g}},u)-\bar{l}_V(y,\bar{l}_V(x,u)) \\
&= \bar{l}_V(x,[s(y),u]_{\wedge})-[s([x,y]_{\mathfrak{g}}),u]_{\wedge}-\bar{l}_V(y,[s(x),u]_{\wedge})\\
&=[s(x),[s(y),u]_{\wedge}]_{\wedge}-[s([x,y]_{\mathfrak{g}}),u]_{\wedge}-[s(y),[s(x),u]_{\wedge}]_{\wedge}\\
&=[s(x),[s(y),u]_{\wedge}]_{\wedge}-[[s(x),s(y)]_{\wedge},u]_{\wedge}-[s(y),[s(x),u]_{\wedge}]_{\wedge}\\
&=0.
\end{align*}
Similarly, we can show that 
\[\bar{l}_V(x,\bar{r}_V(u,y))=\bar{r}_V(\bar{l}_V(x,u),y)+\bar{r}_V(u,[x,y]_{\mathfrak{g}})\]
\[ \bar{r}_V(u,[x,y]_{\mathfrak{g}}= \bar{r}_V(\bar{r}_V(u,x),y)+\bar{l}_V(x,\bar{r}_V(u,y)),\]
for all $x,y \in \mathfrak{g}$ and $u\in V.$ Hence $(V,\bar{l}_V, \bar{r}_{V})$ is a representation of a Leibniz algebra $(\mathfrak{g}, [~,~]_{\mathfrak{g}}).$
Now, $s(T(x))-\hat{T}(s(x)) \in V$ hence $[s(T(x)),u]=[\hat{T}(s(x)),u]$ for all $x\in \mathfrak{g}, u \in V$. Therefore, we have
\begin{align*}
\bar{l}_V&(T(x),T_V(u))=[s(T(x)), T_V(u)]_{\wedge}=[\hat{T}(s(x)),\hat{T}(u)]_{\wedge}\\
&=\hat{T}([\hat{T}(s(x)),u]_{\wedge}+[s(x),\hat{T}(u)]_{\wedge})\\
&=T_V([s(T(x)),u]_{\wedge}+[s(x),T_V(u)]_{\wedge})\\
&=T_V(\bar{l}_V(T(x),u)+\bar{l}_V(x,T_V(u)))
\end{align*}

Hence, $\bar{l}_V(T(x),T_V(u))=T_V(\bar{l}_V(T(x),u)+\bar{l}_V(x,T_V(u)))$ for all $x,y \in \mathfrak{g}$ and $u\in V.$

Similarly, it can be shown that
\[\bar{r}_V(T_V(u),T(x))=T_V(\bar{r}_V(T_V(u),x)+\bar{r}_V(u,T(x))),\]
for all $x \in \mathfrak{g}$ and $u \in V.$ Hence $(V,\bar{l}_V,\bar{r}_V,T_V)$ is a representation of Rota-Baxter Leibniz algebra  $(\mathfrak{g}_T,[~,~]_{\mathfrak{g}}).$
\end{proof}
\begin{proposition}
Let $(\hat{\mathfrak{g}}_{\hat{T}},[~,~]_{\wedge})$ be an abelian extension of $(\mathfrak{g}_T,[~,~]_{\mathfrak{g}})$ by $(V_{T_V},\mu)$. Then any two distinct sections $s_1,s_2: \mathfrak{g} \rightarrow \hat{\mathfrak{g}}$ give the same Rota-Baxter Leibniz algebra representation $(V,\bar{l}_V,\bar{r}_V,T_V).$ 
\end{proposition}
\begin{proof}
Let $s_1$ and $s_2$ be two distinct  sections. Define $\gamma : \mathfrak{g} \rightarrow V$ by
\[\gamma(x)=s_1(x)-s_2(x),~\mbox{for all }~ x \in \mathfrak{g}.\]
Since $\mu (u,v)=0$ for all $ u,v \in V.$
Thus, 
\begin{align*}
[s_1(x),u]_{\wedge}=[\gamma(x)+s_2(x),u]_{\wedge}=[\gamma(x),u]_{\wedge}+[s_2(x),u]_{\wedge}=[s_2(x),u]_{\wedge}
\end{align*}
Similarly, $[u,s_1(x)]_{\wedge}=[u,s_2(x)]_{\wedge}$ for all $x,y \in \mathfrak{g}, u \in V.$
Thus, two distinct sections give the same representation.
\end{proof}
\begin{definition}
Define two maps $\psi : \mathfrak{g}\otimes \mathfrak{g} \rightarrow V$ and $\chi : \mathfrak{g} \rightarrow V$ by  
\[
\psi (x \otimes y)=[s(x),s(y)]_{\wedge}-s([x,y]_{\mathfrak{g}})
\quad \textup{and} \quad
\chi (x)=\hat{T}(s(x))-s(T(x)).
\]
\end{definition}
\begin{proposition}
The cohomological class of $(\psi,\chi)$ does not depend on the choice of sections.
\end{proposition}
\begin{proof}
Let $s_1$ and $s_2$ be two distinct  sections. Define $\gamma : \mathfrak{g} \rightarrow V$ by
\[\gamma(x)=s_1(x)-s_2(x),~\mbox{for all }~ x \in \mathfrak{g}\]
Since $\mu (u,v)=0$ for all $ u,v \in V.$
Therefore, for all $x,y \in \mathfrak{g}, u \in V,$
\begin{align*}
\psi_1(x,y)
&=[s_1(x),s_1(y)]_{\wedge}-s_1([x,y]_{\mathfrak{g}})\\
&=[s_2(x)+\gamma(x),s_2(y)+\gamma(y)]_{\wedge}-(s_2([x,y]_{\mathfrak{g}})+\gamma ([x,y]_{\mathfrak{g}}))\\
&=[s_2(x),s_2(y)]_{\wedge}-s_2([x,y]_{\mathfrak{g}})+[\gamma(x),s_2(y)]_{\wedge}+[s_2(x),\gamma (y)]_{\wedge}-\gamma([x,y]_{\mathfrak{g}})\\
&=\psi_2(x,y)+\delta^1 (\gamma)(x,y).
\end{align*}
Again,
\begin{align*}
\chi_1(x)
&=\hat{T}(s_1(x))-s_1(T(x))\\
&=\hat{T}(s_2(x)+\gamma(x))-s_2(T(x))-\gamma(T(x))\\
&=\chi_2(x)+T_V(\gamma(x))-\gamma (T(x))\\
&=\chi_2(x)-\phi^1(\gamma)(x).
\end{align*}
Therefore, $(\psi_1,\chi_1)-(\psi_2,\chi_2)=(\delta^1 (\gamma),~-\phi^1(\gamma))=d^1(\gamma)$. Hence $(\psi_1,\chi_1)$ and $(\psi_2,\chi_2)$ are in the same cohomology class $H^2_{RBLA}(\mathfrak{g},V)$.
\end{proof}
\begin{theorem}
Let $V$ be a vector space and $T_V:V \rightarrow V$ be a linear map. Then $(V_{T_V},\mu)$ is a Rota-Baxter Leibniz algebra with the bracket $\mu (u,v)=0 $ for all $u,v \in V.$ Then two isomorphic abelian extensions of a Rota-Baxter Leibniz algebra $(\mathfrak{g}_T,[~,~]_{\mathfrak{g}})$ by $(V_{T_V},\mu)$ give rise to the same element in $H^2_{RBLA}(\mathfrak{g},V).$
\end{theorem}
\begin{proof}
Let $(\hat{\mathfrak{g}}_{\hat{T_1}},[~,~]_{\wedge_1})$ and $(\hat{\mathfrak{g}}_{\hat{T_2}},[~,~]_{\wedge_2})$ be two isomorphic abelian extension of $(\mathfrak{g}_T,[~,~]_{\mathfrak{g}})$ by $(V_T,\mu)$. Let $s_1$ be a section of $(\hat{\mathfrak{g}}_{\hat{T_1}},[~,~]_{\wedge_1})$. Thus, we have 
$p_2 \circ (\xi \circ s_1)=p_1 \circ s_1 =Id_{\mathfrak{g}}$ as $p_2 \circ \xi =p_1$, where $\xi$ is the map between the two abelian extensions. Hence $\xi \circ s_1$ is a section of $(\hat{\mathfrak{g}}_{\hat{T_2}},[~,~]_{\wedge_2})$.
Now define $s_2 :=\xi \circ s_1$. Since $\xi$ is a homomorphism of Rota-Baxter Leibniz algebras such that 
$\xi|_{V}=Id_V, ~\xi ([s_1(x),u]_{\wedge_1})=[s_2(x),u]_{\wedge_2}$ . Thus, $\xi|_{V}: V \rightarrow V$ is compatible with the induced representations.

Now, for all $x,y \in \mathfrak{g}$
\begin{align*}
\psi_2(x \otimes y)
&=[s_2(x),s_2(y)]_{\wedge_2}-s_2([x,y]_{\mathfrak{g}})\\
&=[\xi(s_1(x)),\xi(s_1(y)])_{\wedge_2}-\xi(s_1([x,y]_{\mathfrak{g}}))\\
&=\xi ([s_1(x),s_1(y)]_{\wedge_1}-s_1([x,y]_{\mathfrak{g}}))\\
&=\xi (\psi_1 (x \otimes y))=\psi_1(x\otimes y),
\end{align*}
and
\begin{align*}
\chi_2(x)
&= \hat{T_2}(s_2(x))-s_2(T(x))\\
&=\hat{T_2}(\xi(s_1(x)))-\xi(s_1(T(x)))\\
&=\xi (\hat{T_1 }(s_1(x))-s_1(T(x)))\\
&=\xi(\chi_1(x))=\chi_1(x).
\end{align*}
Therefore, two isomorphic abelian extensions give rise to the same element in $H^2_{RBLA}(\mathfrak{g},V).$
\end{proof}

\subsection*{Acknowledgments}
A part of this work has been done when the second author was visiting Centre de Recerca Mathematica (CRM), Barcelona. The second author expresses his gratitude to The Centre International de Mathématiques Pures et Appliquées (CIMPA), France, and CRM, Barcelona for providing their financial support for the visit.


\EditInfo{November 11, 2022}{January 19, 2023}{Ivan Kaygorodov}

\end{paper}